\documentclass[12pt, a4paper, reqno]{amsart}
\usepackage{amsmath,amsfonts,amssymb,amscd,latexsym, amsthm,amsbsy,upref}

\usepackage{enumerate,enumitem,cancel}
\usepackage{bm}
\usepackage[centertags]{amsmath}
\usepackage{amsfonts}
\usepackage{amsthm}
\usepackage{mathrsfs}
\usepackage{mathabx}
\usepackage{xcolor}
\usepackage{graphicx}
\usepackage{comment}
\usepackage{color}
\usepackage{mathrsfs}
\usepackage[utf8]{inputenc}
\usepackage[T1]{fontenc}
\usepackage{tikz}
\usepackage{float}
\usepackage[mathscr]{euscript}
\usepackage{hyperref}
\usepackage{color}

\theoremstyle{definition}
\newtheorem{thm}{Theorem}[section]

\newtheorem{lem}[thm]{Lemma}

\newtheorem{rmk}[thm]{Remark}

\def\hangbox to #1 #2{\vskip3pt\hangindent #1\noindent \hbox to #1{#2}$\!\!$}

\usepackage[marginratio=1:1,height=584pt,width=445pt,tmargin=117pt]{geometry}

\newtheorem{thmA}{Theorem}

\newtheorem{thmx}{Theorem}

\newtheorem*{thm*}{Main Theorem}

\newtheorem*{cor*}{Corollary}

\newtheorem*{prp*}{Proposition}

\newtheorem*{lem*}{Lemma}

\newtheorem*{rmk*}{Remark}

\newtheorem*{qtn}{Question}

\newcommand{\vertiii}[1]{{\left\vert\kern-0.32ex\left\vert\kern-0.32ex\left\vert #1 \right\vert\kern-0.32ex\right\vert\kern-0.32ex\right\vert}}

\newcommand{\Vertiii}[1]{{\left\vert\kern-0.25ex\left\vert\kern-0.25ex\left\vert #1\right\vert\kern-0.25ex\right\vert\kern-0.25ex\right\vert}}

\newcommand{\vertL}[1]{{\left\vert\kern-0.25ex\left\vert\kern-0.25ex\left\vert }}
\newcommand{\vertR}[1]{{\right\vert\kern-0.3ex\right\vert\kern-0.3ex\right\vert}}

\DeclareMathOperator{\spn}{span}

\DeclareMathOperator{\diam}{diam}

\long\def\symbolfootnote[#1]#2{\begingroup%
\def\thefootnote{\fnsymbol{footnote}}\footnote[#1]{#2}\endgroup}

\newcommand{\co}{\mathrm{c}_0}

\allowdisplaybreaks
\begin{document}

\title[Zippin's Embedding Theorem]{Isometric embedding of Banach spaces under optimal projection constants}

\author[C. S. Barroso]{Cleon S. Barroso}
\address[C.S.~ Barroso]{Department of Mathematics, Academic Unit -- Science Center, Federal University of Cear\'a, Fortaleza -- CE, 60455-360, Brazil}
\email{cleonbar@mat.ufc.br}

\date{\today}

\thanks{The research was partially supported by FUNCAP/CNPq/PRONEX Grant 00068.01.00/15.  This work was started when the author was visiting Texas A$\&$M University, Summer 2017. He would like to express his deep gratitude to Prof. Th. Schlumprecht and the Mathematics Department. }


\begin{abstract}
Let $X$ be a Banach space with separable dual. It is proved that for every $\varepsilon\in (0,1)$, $X$ embeds isometrically into a Banach space $W$ with a shrinking basis $(w_n)$ which is $(1+ \varepsilon)$-monotone. Moreover, if $X$ has further an FDD $(E_n)$ whose strong bimonotonicity projection constant is not larger than $\mathcal{D}$, then $(w_n)$ has strong bimonotonicity projection constant not exceeding $\mathcal{D}(1 +\varepsilon)$. Further, if $(E_n)$ is $\mathcal{C}$-unconditional then $(w_n)$ is $\mathcal{C}(1 + \varepsilon)$-unconditional. The proof uses renorming and skipped blocking decomposition techniques. As an application, we prove that every Banach space having a shrinking $\mathcal{D}$-unconditional basis with $\mathcal{D}<\sqrt{6}-1$, has the weak fixed point property. 
\end{abstract}

\maketitle

\setcounter{tocdepth}{1}

\section{Introduction}\label{sec:1}
A classical result of Banach and Mazur \cite[p. 185]{B} states that every separable Banach space embeds isometrically into $C[0,1]$, that is, $C[0,1]$ is isometrically universal for the class of all separable spaces. This remarkable theorem raised an important question in Banach space theory: Does a Banach space satisfying some property embed into a Banach space with a basis having some related and desirable property. In 1988, Zippin \cite{Z2} answering a question originally posed by Pe\l czy\'nski in 1964 \cite[Problem I, p. 133]{Pel}, proved the following fundamental embedding results. 

\begin{thmA}\label{thm:ZI} Every separable reflexive Banach space embeds into a reflexive Banach space with a Schauder basis.  
\end{thmA}

\begin{thmA}\label{thm:ZII} Every Banach space with a separable dual embeds into a space with a shrinking Schauder basis. 
\end{thmA}

Remarkably Theorem \ref{thm:ZII} also solved an important question posed by Lindenstrauss and Tzafriri \cite{LT}. There are nowadays many embedding theorems of that type. Among others, we mention Johnson and Zheng \cite{JoZ}, Odell and Schlumprecht \cite{OS}, Ghoussoub, Maurey and Schachermayer \cite{GMS}, Dodos and Ferenczi \cite{DF}. In \cite{S} Schlumprecht  used coordinate system methods and obtained new proofs for Theorems \ref{thm:ZI} and \ref{thm:ZII}, with more accurate information on the target space. Notably, he proved that every separable dual Banach space $X$ embeds into a space $W$ with a shrinking basis whose Szlenk index of $X$ and the Szlenk index of $W$ are the same. The latest developments in embedding theory suggest quite natural questions. For instance, is it possible to obtain isometric embeddings into spaces with optimal projection constants? Some recent results \cite{Ku, MV, S} hint that at least the monotonicity projection constant can become as optimal as we wish. In 2016, e.g., Kurka \cite{Ku} proved (among other results) that every Banach space with separable dual embeds isometrically into a Banach space with a shrinking monotone basis, and every separable reflexive space embeds isometrically into a reflexive space with a shrinking monotone basis. His argument heavily relies on interpolation methods and techniques from descriptive set theory. The main points can be roughly summarized as follows. Let ${\bf 1}$ be the constant function on the Cantor set $2^\mathbb{N}$ and take $\textrm{g}_0\in C(2^\mathbb{N})$ a function that separates points in $2^\mathbb{N}$. Set $E_X:=\overline{\spn( X \cup \{ {\bf 1}, \textrm{g}_0\})}$. If $X$ is separable and reflexive, then $E_X$ is a separable dual reflexive subspace of the space $C(2^\mathbb{N})$. From descriptive set theory one then gets a shrinking monotone basis $(e^X_i)$ for $E_X$ (cf. \cite[Theorem 5.17]{D}). Finally, applying the interpolation method of Davis, Figiel, Johnson and Pe\l czy\'nski \cite{DFJP}, a new norm on $E_X$ can be built under which $(e^X_i)$ becomes a shrinking and monotone basis, and yet $E_X$ contains an isometric copy of $X$ (\cite[Lemma 4]{Ku}).  

The main goal of this work is to study another related question about embedding theorems. More precisely, we are concerned with the following. 

\begin{qtn} Let $X$ be a Banach space having a Schauder basis $(e_i)$. Does $X$ embed isometrically into a Banach space with a {\it quasi}-monotone basis such that its {\it strong} bimonotonicity constant is of the same magnitude as that of $(e_i)?$ 
\end{qtn}

Recall the strong bimonotonicity projection constant of a Schauder basis $(e_i)$ is the number 
\[
\sup_{m< n}\max\Big( \| P_{[m,n]}\| , \| I - P_{[m,n]}\|\Big),
\]
where for $A\subset \mathbb{N}$, $P_A$ denotes the basis projection onto $\spn(E_i\colon i\in A)$ and $I$ is the identity operator on $X$. Here $[m,n]$ denotes the interval $\{ m, m+1, \dots, n\}$ in $\mathbb{N}$. The basis $(e_i)$ is called $\mathcal{K}$-monotone if $\| P_n\| \leq \mathcal{K}$ for all $n\in \mathbb{N}$ ($P_n:=P_{[1, n]}$). If $\mathcal{K}$ can be made arbitrarily close to $1$, then $(e_i)$ is called {\it quasi}-monotone.

The main result of this paper yields a satisfactory answer for this question. It covers the class of shrinking basis and can be seen as a monotonicity reduction theorem. In particular, the assumption of being shrinking is satisfied when $X$ is reflexive. Although this class seems a little restrictive, as we shall see, it is useful in metric fixed point theory.

%%%%%%%%%%%%%%%%%%%%%%%%%%%%%%%%%%%% SECTION 2

\section{Results}\label{sec:2}
Our main result is the following. 

\begin{thm*}\label{thm:I} For every $\varepsilon>0$, every Banach space $X$ with a separable dual can be isometrically embedded into a Banach space $W$ which has a shrinking monotone basis $(w_n)$ which is $(1+\varepsilon)$-monotone. Moreover:
\begin{itemize}
\item[(i)] If $X$ has an FDD whose strong bimonotonicity constant is not larger than a constant $\mathcal{D}\geq 1$, then the strong bimonotonicity constant of $(w_n)$ is not larger than $\mathcal{D}(1+\varepsilon)$. 
\item[(ii)] If $X$ has an unconditional FDD whose unconditional constant is not larger than a constant $\mathcal{C}\geq 1$, then $(w_n)$ is unconditional with unconditional constant not larger than $\mathcal{C}(1 + \varepsilon)$. 
\end{itemize} 
\end{thm*}

As far as we know, no mention of this result has been made in the literature. Let us point out that FDDs is a natural generalization of Schauder basis. However, in general, it may happen that a Banach space may have a FDD but it may not have a basis (cf. \cite[Theorem 1.1]{Sz}). For the notion of FDD and related terminologies, we refer the reader to Section \ref{sec:3}. 

We will not provide a detailed proof of this result, since it is implicitly contained in the existing literature and can be accomplished through appropriate adaptations. For this reason, we shall only indicate how it follows from \cite{S}. 

\begin{proof}[Proof Sketch]
First, in the proof of \cite[Lemma 2.3]{S} we can simply choose $\rho>1$ close enough to $1$. This leads us to slightly refined versions of the lemmas \cite[Lemma 2.3, Lemma 3.2]{S}. For example, the set $B^*$ in \cite[Lemma 3.2]{S} becomes $\varepsilon$-dense in $B_{X^*}$. These lemmas, when then combined with a skipped blocking procedure, enables us to adapt the construction performed in \cite{S} to embed isometrically $X$  into a Banach space $W$ which has a shrinking basis $(w_n)$ which is $(1+\varepsilon)$-monotone. This yields therefore the first part of the result. It is worth noting, in passing, that the almost-isometric feature of the monotonicity constant of $(w_n)$ is implicitly contained in \cite[Lemmas 2.3 and 2.4]{S} arguments, although it has not (unfortunately) been explicitly observed there. On the other hand, the embeddings constructed in \cite[Proposition 3.3 and Theorem 3.9]{S} are not, however, necessarily isometric. In order to get isometric embeddings, one should thus apply (after having readjusted factors arbitrarily close to $1$, which is necessary) the following well-known lemma (see \cite[Lemma 4.0.10]{Cow}).

\begin{lem*}\label{lem:0} Suppose that $X$ and $Y$ are normed spaces and that $T\colon X\to Y$ is an isomorphic embedding from $X$ into $Y$. Assume that for some $M>0$ one has
\[
\| x\|_X\leq M \| Tx\|_Y\quad\text{ for all } x\in X.
\]
Then the quantity
\[
| y| =\inf\Big\{ \| x\|_X  + M\| y - Tx\|_Y\,\colon \, x\in X\Big\}
\]
defines an equivalent norm on $Y$ with respect to which $T$ is an isometry. Furthermore, for all $y\in Y$, one has
\begin{equation}\label{eqn:1Cowell}
\frac{1}{\|T\|}\|y\|_Y\leq |y|\leq M \| y\|_Y. 
\end{equation}
\end{lem*}

Second, assuming that $X$ has further an FDD with the properties prescribed in assertions (i) and (ii), the superspace $W$ constructed in \cite{S} has the same properties (up to changes by a factor arbitrarily close to 1). The details of this fact can be carried out along the same lines as in \cite[pp. 851--854]{S}. This concludes the sketch of the proof of the theorem. 
\end{proof}

The following result is an immediate consequence of the proof of Main Theorem. 

\begin{thmx}\label{thm:M1} Let $X$ be a Banach space with a shrinking FDD $(E_n)$ and $\mathcal{D}\geq 1$ its strong bimonotonicity constant. Then for every $\varepsilon>0$, $X$ embeds isometrically into a Banach space $W$ with a shrinking Schauder basis $(w_n)$ which is $(1+\varepsilon)$-monotone and has strong bimonotonicity constant at most $\mathcal{D}(1 + \varepsilon)$. 
\end{thmx}

\begin{rmk*} Let us stress that the construction in \cite{S} strongly uses the condition of $1$-norming on a certain coordinate system $(E_j)$ (FDD in our context). Hence, in order to reduce monotonicity and also almost preserve strong bimonotonicity,  the shrinking assumption in Theorem \ref{thm:M1} must be required. Further, it should be noted that every Banach space with an unconditional FDD is isomorphic to a subspace of a Banach space with an unconditional basis (cf. \cite[Theorem 1.g.5]{LT}). 
\end{rmk*}

In the sequel, we will show how this result can be applied in Metric Fixed Point Theory. One of the major concerns in this research area concerns the weak fixed point property (weak-FPP, for short). Let $C$ be a weakly compact convex subset of a Banach space $X$. Does every nonexpansive (i.e., $1$-Lipschitz) mapping $T\colon  C \to C$ have a fixed point$?$ When this happen, we say that C has the weak-FPP. Further, we say that $X$ has the weak-FPP if every weakly compact convex subset of $X$ has the weak-FPP. 

\smallskip 

In 1985, P.-K. Lin \cite{Lin} proved the following beautiful result:

\begin{thmA} Let $X$ be Banach space with a $\mathcal{D}$-unconditional basis satisfying $\mathcal{D}< (\sqrt{33}- 3)/5$. Then $X$ has the weak fixed point property.
\end{thmA}

Since then, significant attention has been paid to  the weak-FPP in Banach spaces with Schauder basis. In particular, it is still an open question as to whether or not every Banach space with an unconditional basis has the weak-FPP. For this setting, the strategy of transferring fixed point analysis to the framework of [X] (the ultrapower of X) has always been a fruitful strategy for studying fixed point property in spaces with Schauder basis. As a consequence, if some quantities related to basis projections can become conveniently small, then Lin's approach \cite{Lin} can be successfully applied to obtain better results. This is the key idea behind our second result (Theorem \ref{thm:M2}). Note that the first statement in our main result implies that the bimonotonicity constant of $W$ can be kept under $2+\varepsilon$, consequently the strong bimonotonicity constant can be kept under $3+\varepsilon$. However, these quantities are not small enough for applications in fixed point theory. Nevertheless, Theorem \ref{thm:M1}  allows us to obtain an interesting improvement of Lin's result, since the monotonicity constant can be as close as possible to $1$ and the strong bimonotonicity constant is almost preserved. 

\begin{thmx}\label{thm:M2} Let $X$ be a Banach space with a shrinking unconditional basis whose unconditional constant satisfies $\mathcal{D}< \sqrt{6}-1$. Then $X$ has the weak fixed point property. 
\end{thmx}

Theorem \ref{thm:M2} immediately provides the following.

\begin{thmx}\label{thm:M3} Every reflexive Banach space with a $\mathcal{D}$-unconditional basis satisfying $\mathcal{D}< \sqrt{6}-1$ has the fixed point property. 
\end{thmx}

In general, a Banach space $X$ is said to have the fixed point property when every bounded, closed convex subset $C\subset X$ has the property that every nonexpansive self-mapping of $C$ has a fixed point. In light of the above results, one can speculate whether Theorem \ref{thm:M2} remains valid without the assumption of shrinking. This seems however to require another approach as the proof of Theorem \ref{thm:M1} uses this condition to build $1$-norming FDDs through skipped blocking decompositions. It worths to point out that an FDD of a Banach space with a separable dual needs not be shrinking. In fact, a result of Zippin \cite{Z1} ensures, e.g., that $\co$ has a nonshrinking basis (a fortiori a nonshrinking FDD), since $\co$ is nonreflexive. Despite that, it can be replaced by a more weaker one, as our last result shows. 

\begin{thmx}\label{thm:M4} Every Banach space with a $\mathcal{D}$-unconditional spreading basis satisfying $\mathcal{D}<\sqrt{6}-1$ has the weak fixed point property. 
\end{thmx}

\smallskip  

This paper is organized as follows. In Section \ref{sec:3}, we recall some basic definitions and set up the main ingredients used throughout the paper. In Section \ref{sec:4} we present the proofs of Theorems \ref{thm:M2} and \ref{thm:M4}.

%%%%%%%%%%%%%%%%%%%%%%%%%%%%%%%% SECTION 3

%%%%%. PRELIMINARIES

%%%%%%%%%%%%%%%%%%%%%%%%%%%%%%%%%%%%%%

\section{Background}\label{sec:3}
We will use standard Banach space notation and terminology, mostly contained e.g. in \cite{AK, HJ, S}.  Let $X$ be a Banach space. A sequence $(e_n)_{n=1}^\infty$ in $X$ is said to be a basis of $X$ is for each $x\in X$ there is a unique sequence of scalars $(a_n)_{n=1}^\infty$ such that $x= \sum_{n=1}^\infty a_n e_n$, that is, $\| x - \sum_{i=1}^n a_i e_i\| \to 0$. A sequence $(e_n)_{n=1}^\infty \subset X$ is called a {\it basic sequence} if it is a basis for its closed linear span $[e_n\colon n\in \mathbb{N}]$. A series $\sum_{n=1}^\infty x_n$ in $X$ is said to {\it converge unconditionally} precisely when $\sum_{n=1}^\infty x_{\pi(n)}$ converges for every permutation $\pi$ of $\mathbb{N}$. A Schauder basis $(e_i)_{i=1}^\infty$ of a Banach space $X$ is said to be {\it unconditional} if for every $x\in X$, its expansion $x=\sum_{i=1}^\infty a_i e_i$ converges unconditionally. For more on unconditional bases, we refer the reader to \cite[p. 200]{FHHMZ}

\subsection{Finite Dimensional Decompositions}\label{sec:3.1}
Recall from \cite{JL} that a sequence of finite dimensional subspaces $(Z_j)$ of a Banach space $(Z,\| \cdot\|)$ is called a finite dimensional decomposition (FDD) of $Z$ if each $z\in Z$ can be uniquely written as $z=\sum_{n=1}^\infty z_n$ with $z_n\in Z_n$ for all $n\in \mathbb{N}$. For $n\in \mathbb{N}$, the $n$-th projection $P^Z_n$ on $Z$ associated to an FDD $(Z_j)$ is defined by $P^Z_n(z)= z_n$. It is well known that the sequence of projections $(P_n^Z)$ satisfies $\sup_n \| P^Z_n\|<\infty$. For a finite set $A\subset \mathbb{N}$, define $P^Z_A=\sum_{n\in A} P^Z_n$. The projection constant of $(Z_j)$ is defined by $\mathcal{K}^Z_b((Z_j))=\sup_{m\leq n} \| P^Z_{[m,n]}\|$. We will also refer to it as {\it bimonotonicity projection constant}. An FDD $(Z_j)$ of $Z$ is called $\mathcal{D}$-{\it bimonotone} $(\mathcal{D}\geq 1)$ if $\mathcal{K}^Z_b( (Z_j))\leq \mathcal{D}$. Recall also that an FDD $(Z_n)$ is called $\mathcal{D}$-monotone if $\mathcal{K}^Z_m((Z_j))=\sup_{n\geq 1}\big\| P^Z_{[1, n]}\big\|\leq \mathcal{D}$. Similarly to Schauder basis, recall that the strong bimonotonicity projection constant of an FDD $(Z_n)$ is the number 
\[
\mathcal{K}^Z_{sb}((Z_j)):=\sup_{m< n}\max\Big( \| P^Z_{[m,n]}\| , \| I - P^Z_{[m,n]}\|\Big),
\]
where for $A\subset \mathbb{N}$, $P^Z_A$ denotes the basis projection onto $\spn(Z_i\colon i\in A)$ and $I$ is the identity operator on $X$. We say that $(Z_j)$ is $\mathcal{D}$-{\it strongly bimonotone} if 
\[
\mathcal{K}^Z_{sb}((Z_j))\leq \mathcal{D}.
\]  
In the special case when $\mathcal{D}=1$, we simply say that $(Z_j)$ is {\it monotone}, {\it bimonotone} or {\it strongly bimonotone}, respectively.

Note that if an FDD $(Z_j)$  is $\mathcal{D}$-{\it strongly bimonotone} then
\[
\Bigg\| \sum_{j=1}^m x_j + \sum_{j=n}^\infty x_j\Bigg\| \leq \mathcal{D}\Bigg\| \sum_{j=1}^\infty x_j\Bigg\|.
\] 
for all $x=\sum_{j=1}^\infty x_j \in \spn(Z_j\colon j\in \mathbb{N})$ and for all $m<n$ in $\mathbb{N}$.  An FDD $(Z_j)$ is called {\it unconditional} if $z= \sum_{n=1}^\infty z_n$ converges unconditionally, for all $z\in Z$. In this case, similarly to basis, given a sequence of signs $(\epsilon_n)_{n=1}^\infty \in \{ -1, 1\}^{\mathbb{N}}$, the series $\sum_{n=1}^\infty \epsilon_n z_n$ converges. Further, the Uniform Boundedness Principle implies that 
\[
\mathcal{C}_u= \sup\Bigg\{ \Bigg\| \sum_{n=1}^\infty \epsilon_n z_n\Bigg\| \colon (\epsilon_n)_{n=1}^\infty \in \{ -1, 1\}^{\mathbb{N}}\Bigg\}<\infty
\]
The quantity $\mathcal{C}_u$ is called the {\it unconditional constant} of the FDD.  

\begin{rmk}\label{rmk:3.2} Clearly, if  an FDD $(Z_j)$ of $Z$ is $\mathcal{D}$-strongly bimonotone then its projection constant satisfies $\mathcal{K}_b((E_j))\leq \mathcal{D}$. It is known (and easy to see) that $\mathcal{D}$-unconditional bases are suppression $(\mathcal{D}+1)/2$-unconditional. A proof of this fact is implicitly contained in the proof of \cite[(ii)$\to$(iii), p.202]{FHHMZ}. It follows therefore that every shrinking $\mathcal{D}$-unconditional basis gives rise to a shrinking $(\mathcal{D}+1)/2$-strongly bimonotone FDD. 
\end{rmk}

\subsection{Preliminaries on Metric Fixed Point Theory} Before going to the proofs of Theorems \ref{thm:M2} and \ref{thm:M4} we fix some notation. Let $\mathscr{U}$ stand for a free ultrafilter on $\mathbb{N}$. Denote by $\lim_{i\to\mathscr{U}}$ the ultra-limit through $\mathscr{U}$. For a Banach space $X$, the $\ell_\infty$-sum of $X$ is given by $\ell_\infty(X)=\{ (v_i) \colon \text{each } v_i\in X\,\text{ and }\, \sup_{i\in \mathbb{N}}\| v_i\|<\infty\}$. The ultrapower $[X]$ of $X$ is the quotient space $\ell_\infty(X)/ \mathcal{N}$ equipped with the norm $\| [v_i]\|=\lim_{i\to \mathscr{U}}\| v_i\|$, where $\mathcal{N}= \{ (v_i)\in \ell_\infty(X)\colon \lim_{i\to \mathscr{U}}\|v_i\|=0\}$ and $[v_i]$ denotes the element of $[X]$ associated to the sequence $(v_i)$. Clearly, $X$ embeds isometrically into $[X]$ through the map $x\mapsto [x, x,\dots]\in [X]$, so henceforth we will not distinguish between $x$ and $[x,x,\dots]$. If $C$ is a weakly compact convex subset of $X$ and $T\colon C\to C$ is a nonexpansive fixed-point free mapping, then by Zorn's lemma there is a weakly compact convex set $K\subset C$, $\diam K>0$, which is minimally invariant under $T$. In this context, Goebel \cite{G} and Karlovitz \cite{Ka} proved the following result.

\begin{lem} Let $T\colon K\to K$ be a nonexpanisve fixed-point free mapping on a minimal weakly compact convex set. Assume that $(y_n)$ is an approximate fixed point sequence of $T$. Then 
\[
\lim_{n\to \infty} \| x - y_n\| =\diam K\quad\text{ for all } x\in K. 
\]
\end{lem}
Recall that a sequence $(y_n)\subset K$ is called an approximate fixed point sequence of $T$ provided that $\| y_n - Ty_n\| \to 0$. In the ultra-product language, both $K$ and $T$ are described as follows 
\[
[K]=\big\{ [v_i] \in [X] \colon v_i \in K\, \forall\, i\in \mathbb{N}\big\}
\]
and $[T]( [v_i]):= [Tv_i]$ for all $[v_i]\in [K]$. The mapping $[T]$ is called the ultrapower mapping induced by $T$. It is easy to see that $[T]$ is nonexpansive and maps $[K]$ into itself. Moreover, the set of fixed points of $[T]$ is nonempty and consists of all points $[x_i]\in [ K]$ for which $(x_i)$ is an approximate fixed point sequence of $T$. The ultrapower counterpart of Goebel-karlovitz's lemma (see also \cite[Corollary 3.2]{AkK}) is the following result due to P.-K. Lin \cite{Lin}. 

\begin{lem}\label{lem:Lin} Let $K$ be a minimal weakly compact convex set for a nonexpansive fixed-point free mapping $T$. Assume that $[M]$ is a convex nonempty subset of $[K]$ which is invariant under $[T]$. Then for all $x\in K$, 
\[
\sup\{ \| x - [v_i] \| \,\colon \, [v_i] \in [M]\} = \diam K. 
\]
\end{lem}

%%%%%%%%%%%%%%%%%%%%%%%%%%%%%
\smallskip

%%%%%%%%%%%%%%%%%%%%%%%%%%

\section{Proofs of Theorem \ref{thm:M2} and \ref{thm:M4}}\label{sec:4}

\subsection{Proof of Theorem \ref{thm:M2}}
Let $X$ be as in the statement of the theorem and let $(e_i)$ denote its $\mathcal{D}$-unconditional basis in which $\mathcal{D}< \sqrt{6}-1$. Towards a contradiction, assume that $X$ fails the weak fixed point property. Then there is a weakly compact convex subset $K$ of $X$ which is minimally invariant under a nonexpansive fixed-point free mapping $T$. Fix $(y_n)$ an approximate fixed point sequence of $T$. By Goebel-Karlovitz's lemma, we have 
\[
\lim_{n\to \infty} \| x - y_n\|=\diam K\qquad \text{for all } x\in K. 
\]
By passing to a subsequence and considering a translation followed by an appropriate scaling of $K$, we may assume that $\diam K=1$ and $(y_n)$ is weakly convergent to $0$. So, $(y_n)$ is a seminormalized weakly null sequence. Fix $\varepsilon>0$ small enough so that 
\begin{equation}\label{eqn:1sec10}
\frac{\varepsilon}{2}\cdot \Big[ \frac{ (\mathcal{D} +1)( \mathcal{D} + 2 + \varepsilon)}{2} +1 \Big] < \frac{3}{2} - \frac{ (\mathcal{D} +1)^2 }{4}. 
\end{equation}
This is possible because $\mathcal{D}< \sqrt{6}-1$. Since $(e_i)$ is $\mathcal{D}$-unconditional, every block basis of $(e_i)$ is $\mathcal{D}$-unconditional, too. Thus, using a standard gliding hump argument, we can select a subsequence $(x_n)$ of $(y_n)$ which is $\mathcal{D} + \varepsilon$-unconditional. It follows that $\| x_m + x_n\| \leq (\mathcal{D} + \varepsilon) \| x_m - x_n\|$ for all $m< n$ in $\mathbb{N}$. As $\diam K=1$, we get
\begin{equation}\label{eqn:2sec10}
\| x_m + x_n\| \leq \mathcal{D} + \varepsilon\quad \text{for all } m< n.
\end{equation}
Recall from Remark \ref{rmk:3.2} that $(e_i)$ is suppression $(\mathcal{D}+1)/2$-unconditional. In particular, it is $(\mathcal{D}+1)/2$-strongly bimonotone. By the Main Theorem, there is an isometric linear embedding $J$ from $X$ into a Banach space $W$ with a shrinking $(1+\varepsilon)$-monotone basis $(w_i)$ which is $(\mathcal{D} +1)(1 +\varepsilon)/2$-strongly bimonotone. Denote by $\vertiii{\cdot }$ the norm of $W$. It follows from (\ref{eqn:2sec10}) that 
\begin{equation}\label{eqn:3sec10}
\vertiii{ J x_m + Jx_n} \leq \mathcal{D} +\varepsilon\quad \text{for all } m< n.
\end{equation}
Set $K^J=J(K)$ and define $T^J\colon K^J\to K^J$ by $T^J( Jx)= J(Tx)$. Then $K^J$ is a minimal weakly compact convex set for the nonexpansive fixed-point free mapping $T^J$. Further, one also has $\diam K^J=1$. The next step of the proof concerns a standard procedure in metric fixed point theory. Let $(P_i)$ denote the natural projections associated to the basis $(w_i)$. For $n\in \mathbb{N}$, set $R_n= I- P_n$. Let $[W]$ denote the ultrapower of $W$. Since $(Jx_n)$ is weakly null, we can use the gliding hump argument to find a subsequence $(x_{m_i})$ of $(x_n)$ and a sequence of non-overlapping intervals $(F_i)$ in $[\mathbb{N}]^{<\infty}$ satisfying the following properties:
\begin{itemize}
\item[($P1$)] $\vertiii{ [J x_{m_i} ] - [ P_{F_i} Jx_{m_i} ]}=0$.\vskip .15cm 
\item[($P2$)] $\vertiii{ [ Jx_{m_{i+2}}] - [ R_{\max F_i} Jx_{m_{i+2}} ] } =0$. \vskip .15cm 
\item[$(P3$)] $\vertiii{ [ P_{F_i} Jx ]}= \vertiii{ [ R_{\max F_i} Jx ] } =0$ for all $x\in X$. 
\end{itemize}
Set $\tilde{x}= [ Jx_{m_i}]$ and $\tilde{y}= [ Jx_{m_{i+2}}]$. Now consider the set
\begin{align*}
[M]=\left\{  [v_i]\in [K^J]\colon\, \begin{matrix}&\exists\, u\in K\text{ such that }\, \vertiii{ [v_i] - Ju}\leq (\mathcal{D}+\varepsilon)/2, \text{ and }\\ 
&\max\big( \vertiii{ [ v_i] - \tilde{x}} , \vertiii{ [v_i] - \tilde{y}}\big)\leq 1/2
\end{matrix}\right\}.
\end{align*}
Let $[T^J]\colon [K^J]\to [K^J]$ be the ultrapower mapping induced by $T^J$. Direct calculation shows that $[T^J]$ is nonexpansive and leaves $[M]$ invariant. Then Lemma \ref{lem:Lin} implies
\begin{equation}\label{eqn:4sec10}
\sup\big\{ \vertiii{ [v_i]} \colon \, [v_i]\in [M]\big\}=1. 
\end{equation}
As it turns out, however, this cannot be true. Indeed, first, since $\diam K^J=1$, we deduce from (\ref{eqn:3sec10}) that $(\tilde{x} +\tilde{y})/2\in [M]$. Let now $\tilde{P}=[ P_{F_i}]$ and $\tilde{Q}=[ R_{\max F_i}]$. Our choice for the operator $\tilde{Q}$ was inspired by the work of Khamsi \cite[p. 29]{K}. Fix any $[v_i]\in [M]$ and choose $u\in K$ so that $\vertiii{ [ v_i ] - Ju}\leq (\mathcal{D} + \varepsilon)/2$. Then 
\[
\begin{split}
\vertiii{  [ v_i]} \leq \frac{1}{2}\Bigg( \Vertiii{ \tilde{P} + \tilde{Q}} \vertiii{ [v_i] - Ju} + \vertiii{ [I] - \tilde{P} } \vertiii{ [v_i] - \tilde{x}} + \vertiii{ [ I] - \tilde{Q}} \vertiii{ [v_i] - \tilde{y} } \Bigg).
\end{split}
\]
Using now that $(w_i)$ is $(1+\varepsilon)$-monotone and is also $(\mathcal{D}+1)(1 + \varepsilon)/2$-strongly bimonotone, we conclude that 
\[
\max\Big(\vertiii{ \tilde{P} + \tilde{Q}}, \vertiii{ [I] - \tilde{P} }\Big) \leq \frac{(\mathcal{D}+1)(1+\varepsilon)}{2}\,\,\text{ and }\,\,  \vertiii{ [I] - \tilde{Q} }\leq 1+\varepsilon
\]
It follows that 
\[
\vertiii{  [ v_i]} \leq \frac{1}{2}\Bigg(  \frac{(\mathcal{D}+1)(1+\varepsilon)}{2}\cdot\frac{(\mathcal{D} + \varepsilon)}{2} +  \frac{(\mathcal{D}+1)(1+\varepsilon)}{2}\cdot \frac{1}{2} + \frac{1 + \varepsilon}{2}\Bigg),
\]
and hence, taking the supremum over all  $ [ v_i]\in [M]$, we obtain 
\[
\sup\big\{ \vertiii{ [v_i]} \colon \, [v_i]\in [M]\big\}\leq \frac{1}{2}\Bigg(  \frac{(\mathcal{D}+1)(1+\varepsilon)}{2}\cdot\frac{(\mathcal{D} + \varepsilon)}{2} +  \frac{(\mathcal{D}+1)(1+\varepsilon)}{2}\cdot \frac{1}{2} + \frac{1 + \varepsilon}{2}\Bigg).
\] 
From our choice of $\varepsilon$ in (\ref{eqn:1sec10}) we deduce that $\sup\big\{ \vertiii{ [v_i]} \colon \, [v_i]\in [M]\big\}<1$, contradicting (\ref{eqn:4sec10}). This therefore completes the proof of theorem.   \hfill $\qed$

\medskip 

%%%%%%%%%%%%%%%%%%%%%

\subsection{Proof of Theorem \ref{thm:M4}} Let $X$ be a Banach space with $\mathcal{D}$-unconditional spreading basis $(e_i)$ satisfying $\mathcal{D}< \sqrt{6}-1$. Assume by way of contradiction that $X$ fails the weak fixed point property. In particular, $X$ does not have Schur's property. Since $(e_i)$ is spreading, it follows that $(e_i)$ is a seminormalized weakly null sequence in $X$ (see \cite[p. 199, Fact 37]{HJ}). Fix $\varepsilon>0$ as in (\ref{eqn:1sec10}). By the well-known Bessaga-Pe\l czy\'nski's selection principle, we can find a subsequence $(e_{n_i})$ of $(e_i)$ which is $1+\varepsilon$-basic. Denote $(P_n)$ the natural projections associated with $(e_i)$. Since $(e_i)$ is spreading, we deduce that $\sup_{n\geq 1}\| P_n\|\leq 1+\varepsilon$. Therefore, we can proceed in the same way as in the proof of Theorem \ref{thm:M2} to arrive at a contradiction after applying Lin's Lemma \ref{lem:Lin}. The proof is complete.  \hfill $\qed$

\begin{rmk} Recall that a basis $(e_i)$ is said to be spreading if it is $1$-equivalent to all of its subsequences. 
\end{rmk}

\medskip 
\subsection*{Acknowledgements}\label{ackref} The author would like to thank the anonymous reviewer for very careful proofreading, helpful comments and suggestions leading to a much improved version of the manuscript.

\end{document}